\documentclass[11pt,reqno]{amsart}
\usepackage{graphicx,amsmath,amsfonts,latexsym,amssymb,amsthm,color}
\usepackage{tikz-cd}
\usepackage{verbatim}
\usepackage{enumerate}

\usepackage[a4paper]{geometry}
\definecolor{blau}{rgb}{0.05,0.2,0.7}
\definecolor{auchblau}{rgb}{0.03,0.3,0.7}
\usepackage[pdftex,linkbordercolor={0.8 0.5 0.2}]{hyperref}

\renewcommand{\Re}{\operatorname{Re}}

\newcommand{\der}{\mathrm{d}}

\newcommand{\rmi}{\mathrm{i}}

\newcommand{\R}{\mathbb{R}}
\newcommand{\C}{\mathbb{C}}
\newcommand{\N}{\mathbb{N}}
\newcommand{\calK}{\mathcal{K}}

\newcommand{\calU}{\mathcal{U}}

\newtheorem{theorem}{Theorem}[section]

\newtheorem{example}[theorem]{Example}

\newtheorem{rem}[theorem]{Remark}

\title[Heat kernel expansion]{Heat and Wave kernel expansions for stationary spacetimes}

\author[A.~Strohmaier]{Alexander Strohmaier}
\address{Leibniz University Hannover, Institute of Analysis, 30167 Hannover, Germany}  \email{a.strohmaier@math.uni-hannover.de} 

\author{Steve Zelditch (\textdagger)}

\address{Department of Mathematics, Northwestern  University,
Evanston, IL 60208-2370, USA} \email{
zelditch@math.northwestern.edu}

\thanks{Research partially supported by NSF grant   DMS-1810747
and by the Stefan Bergman trust.}

\begin{document}

\begin{abstract}
The generator of time-translations on the solution space of the wave equation on stationary spacetimes specialises to the square root of the Laplacian on Riemannian manifolds when the spacetime is ultrastatic. Its spectral analysis therefore constitutes a generalization of classical spectral geometry.
If the spacetime is spatially compact the spectrum is discrete and admits a wave-trace expansion at time zero. 
A Weyl law for the eigenvalues and a wave-trace formula was shown in \cite{SZ18} and related to the geometry of the space of null-geodesics. In this paper we investigate the relation to heat kernel coefficients and residues of zeta functions in this context and compute the second non-zero term in the wave-trace expansion. This second coefficient is an analogue in the category of stationary spacetimes of the second heat kernel coefficient of the Laplace operator. The general formula is quite involved but reduces to the usual term involving the scalar curvature when specialised to ultra-static spacetimes.
 \end{abstract}


\maketitle


\section{Introduction and Setting}

Hyperbolic evolution equations on globally hyperbolic spacetimes, such as the Klein-Gordon equation or the wave-equation admit a time-independent Hamiltonian in case there exists a complete timelike Killing field, i.e. in case the spacetime is stationary (see for example \cite{Wal84}). This Hamiltonian can be viewed as an operator and its spectral theory naturally generalizes spectral theory of the Laplacian on Riemannian manifolds. This not only is interesting for purely mathematical reasons; methods of spectral theory and scattering theory have been found to be useful, for example, in general relativity (see for example \cite{MR2426141,MR2864787,MR2215614,MR4213773,MR4365146,MR3038715} and references therein). 
Ultimately spectral properties describe processes of dispersion, decay, and oscillations of physical quantum particles on curved backgrounds.
Typically a Cauchy hypersurface is chosen and separation of variables results in some kind of spectral problem on the space of functions on the Cauchy surface. As argued in previous work (\cite{SZ18,SZ20}) by the authors there is a completely natural formulation of spectral theory for the Klein-Gordon equation on stationary spacetimes, that is invariantly defined and therefore independent of the choice of a Cauchy surface. In this context we proved a Weyl law and a Gutzwiller-Duistermaat-Guillemin-type (\cite{Gutz71,DG75}) trace formula, and we computed the principal wave-trace invariants in terms of invariantly defined quantities of the stationary space-time.

Global hyperbolicity of the spacetime, which is a very natural condition from the point of view of evolution equations, has the topological implication that the spacetime $M$ is diffeomorphic to a product $\R \times \Sigma$ with each $\Sigma_t = \{t\} \times \Sigma$ being a smooth spacelike Cauchy hypersurface. Because of this structure the existence of the complete time-like Killing field poses no additional topological constraint.
To describe the local geometry of stationary spacetimes we give two descriptions which emphasise different geometric aspects. 

The first description depends on the choice of a smooth spacelike Cauchy surface $(\Sigma,h)$ with induced metric $h$. Each Killing orbit, being a timelike curve, must intersect this Cauchy surface exactly once. This defines a diffeomorphism $\Phi: \R \times \Sigma \to M, (t, y) \mapsto \mathrm{exp}(t Z) y$, where $\mathrm{exp}(t Z)$ is the exponential map associated to the Killing field $Z$. 
The pull-back of the metric $g$ with respect to this diffeomorphism takes the form
$$
 \Phi^* g = - N^2 \der t^2 + h_{jk} (\der x ^j + w^j \der t)(\der x ^k + w^k \der t),
$$
with respect to any local chart on $\Sigma$. 
One major difference to the product case is that the future pointing unit normal $\nu$ to the Cauchy surface $\Sigma$ is not $\partial_t$ but
$\partial_\nu = \frac{1}{N} \left(\partial_t - w^j \partial_j \right)$. 
In this expression
$N \in C^\infty(\Sigma)$ is a smooth positive function, $w = w^k \partial_k$ is a vector field on $\Sigma$, and we have used the Einstein sum convention.
The function $N$ and the vector field $w$ are sometimes referred to as the Lapse function and the shift vector field. The product case corresponding to Riemannian spectral geometry is recovered in the special case $N=1, w=0$.
The above formula also provides a large class of examples. Given a  Riemannian manifold $(\Sigma,h)$ with a positive smooth function $N \in C^\infty(\Sigma)$ and a vector field $w$ such that the above expression has Lorentzian signature this defines a Lorentzian metric on $\R \times \Sigma$ with complete Killing field $Z = \partial_t$. 

A second more invariant description is to think of $(M,g)$ as the total space of an $\R$-principal bundle over the space $\calK$ of Killing orbits. By the above description the space of Killing orbits inherits the structure of a smooth manifold from any Cauchy surface, and this structure is independent of the choice of Cauchy surface.
We then have the natural projection map $\pi: M \to \calK$ and the metric takes the form
$$
 g = - u^{2} \theta \otimes \theta + \pi^* g_\calK.
$$
Here $\theta$ is a one-form with $\mathcal{L}_Z \theta =0$, $\theta(Z) =1$, $\theta|_{Z^\perp}=0$, and $u: \calK \to \R$ is a positive smooth function. These properties uniquely determine $\theta$ and the 
Riemannian metric $g_\calK$ on the quotient space such that the horizontal lift of vectors is an isometry. This is the context of (pseudo-) Riemannian principal bundles with projections that are (pseudo-) Riemannian submersions.
The above description in the context of stationary spacetimes has also appeared in the mathematical literature and is related to Finsler geometry; see, for example, \cite{MR2281195}.  We refer to
 \cite{MR2505228} for various different representations of the metric, their relation, and physical interpretations.
 
As explained above any choice of Cauchy surface $\Sigma$ leads to an identification of $\Sigma$ with the space of Killing orbits $\calK$ and the corresponding functions can therefore be related. We have
\begin{align*}
 u^2  = -g(Z,Z) = N^2 - |w|_h^2, \quad u^2 | \theta |_g ^2 =-1,\\
 \theta = \der t + \gamma,  \textrm{ where}\quad \gamma = - u^{-2} h_{jk} w^k \der x^j,
\end{align*}
and $g_\calK$ is identified with $h + u^{-2} w^\flat \otimes w^\flat$, where the musical isomorphism is with respect to the metric $h$. We also note that the quotient admits a natural volume form $\der \mathrm{Vol}_\calK$ defined by contraction. It is explicitly given by
$$
 \der \mathrm{Vol}_\calK (X_1,\ldots,X_{n-1}) = \der \mathrm{Vol}_g(Z,\tilde X_1,\ldots, \tilde X_{n-1}),
$$
where $\tilde X_k$ denotes the horizontal lift of the vector field $X_k$ on $\calK$ defined by $\theta(\tilde X_k)=0$ and $\pi_*(\tilde X_k)= X_k$. This volume form is related to the metric volume form $\der \mathrm{Vol}_{g_\calK}$ on $(\calK,g_\calK)$ by
$\der \mathrm{Vol}_{g_\calK} = u^{-1}  \der \mathrm{Vol}_\calK$.
In terms of the  metric volume form on $\Sigma$ we have
$$
 \der \mathrm{Vol}_\calK = N \cdot \der \mathrm{Vol}_{h}.
$$
The quantities $N \der \mathrm{Vol}_{h}$ and $|Z| := u= (-g(Z,Z))^\frac{1}{2} = (N^2 - |w|_h^2)^{\frac{1}{2}}$ are therefore invariantly defined in terms of the geometry of the stationary spacetime as an $\R$-principal bundle. In contrast the integral of a function on $\calK$ with respect to the volume form $\der \mathrm{Vol}_{h}$  on $\Sigma$ depends on the choice of $\Sigma$. It is not invariant under local gauge transformations  (see Section \ref{invariance}). 

Locally this type of stationary geometry appears in general relativity for example for rotating stationary stars. In the absence of rotation one has $w=0$ and $N$ describes the effects of gravitation. A non-zero shift vector field $w$ adds additional geometry that describes a frame-dragging effect. The Kerr spacetime is not stationary in the sense discussed here due to the existence of the ergo-sphere in which the timelike Killing field fails to be timelike.
However, the far away region in Kerr spacetime is locally stationary and of the standard form. To understand the geometry better and illustrate the theorems we consider the example of a rotating sphere and we choose for simplicity $N=1$.

\begin{example}[A rotating round sphere with $N=1$] \label{examplerot}
Given a number $\Omega \in \R$ with  $|\Omega|<1$ we consider the sphere 
$(\Sigma,h) = (S^2, \der \vartheta^2 + \sin^2(\vartheta)\, \der \varphi^2)$ with its round metric and we choose
$$
 w = \Omega \partial_\varphi.
$$
With $N=1$ the stationary metric on $M=\R \times S^2$ takes the form
$$
 g = - \der t^2 + \der \vartheta^2 + \sin^2(\vartheta)\,(\der \varphi + \Omega \der t)^2.
$$
A complete Killing field on $M$ is $Z=\partial_t$, and
$$
 |Z|^2 = -g(Z,Z)= 1- |w|_h^2 = 1-\Omega^2 \sin^2(\vartheta),
$$
so $Z$ is timelike everywhere precisely when $|\Omega|<1$.
A direct computation shows that
$$
 d\theta
 =
 -\frac{2\Omega\sin(\vartheta)\cos(\vartheta)}{(1-\Omega^2\sin^2(\vartheta))^2}\,
 d\vartheta\wedge d\varphi,
$$
so $\theta$ is closed if and only if $\Omega=0$. Thus for $\Omega\neq 0$ the spacetime is stationary but not static.
This example is particularly simple because the shift vector field $w$ is Killing.
\end{example}

In this paper we focus on wave-trace invariants near time zero which are invariants that in classical spectral geometry give rise to the Weyl law, heat kernel coefficients and the poles of the zeta function. 
As in previous work we will specialize to the case of stationary spacetimes with {\sl compact Cauchy hypersurface}. This class plays the same role as closed Riemannian manifolds do in ordinary spectral theory: they form the building blocks of the theory and provide the foundation to more general treatments via scattering theory. In classical spectral geometry not only  Einstein manifolds are considered. We similarly do not assume our spacetimes to satisfy the vacuum Einstein equations. We would like to point out that spectral geometry on non-compact manifolds is a rich theory in which local computations and theorems from spectral geometry on compact manifolds enter (see \cite{MR3729409} and references for a recent survey). Good examples are Weyl laws for resonances or the spectral shift function which both are reduced to their compact counterparts (see for example \cite{MR0278108, MR0481639, MR0512084, MR0575734,MR1890995}). Such theorems hold in greater generality and we expect the same to be true for the spectral geometry on stationary spacetimes. 

We now describe the main result of the paper. 
We fix a globally hyperbolic stationary and spatially compact spacetime $(M,g)$ of dimension $n \geq 3$ with a complete timelike Killing vector field $Z$. In particular there exists a compact spacelike Cauchy surface.
For a smooth invariant potential $W \in C^\infty(M,\R)$, $\mathcal{L}_Z W=0$ we consider the wave operator $\Box = \Box_g + W(x)$. A typical example is the minimally coupled Klein-Gordon operator
$\Box = \Box_g + m^2 + \kappa\, \mathrm{scal}$, where $\mathrm{scal}$ denotes the scalar curvature and $\kappa \in \R$. 
By invariance the wave operator $\Box$ commutes with the Killing field, and therefore with the pull-back of the Killing flow $e^{s Z}$
acting on functions $C^\infty(M)$. We let $\mathcal{H}^\infty = \mathrm{ker}(\Box) \cap C^\infty(M)$ be the space of smooth solutions of $\Box u =0$ on $M$. 

As proved in \cite{SZ18} there exists a decomposition $\mathcal{H}^\infty = \mathcal{V} \oplus \mathcal{W}$ into $\mathcal{L}_Z$-invariant subspaces such that $\mathcal{V} \subset \mathcal{H}^\infty$ is finite dimensional. Moreover, the infinite dimensional space $\mathcal{W}$ admits a basis (a linearly independent total set), $(u_j)_{j \in \mathbb{Z}}$ consisting of eigenvectors for $\rmi \mathcal{L}_Z$ with real
 and non-zero eigenvalues $\lambda_j$ such that $\lambda_j \leq \lambda_{j+1}$ and such that $\lambda_j \to \pm \infty$ as $j \to \pm \infty$.
In other words the generator $H = \rmi \mathcal{L}_Z$ of the Killing field can be written as $H=H_0 + H_1$, where $H_0$ is an endomorphism of a finite dimensional space and $H_1$ is diagonalizable with eigenvalues $\lambda_j$. Of course $H_0$ can then be further decomposed into invariant root subspaces consisting of generalized eigenvectors. We denote these eigenvalues by $(\beta_k)_{k=1,\ldots,N}$ and we denote by $m_k$ the algebraic multiplicity of $\beta_k$. 
Independent of the choice of splitting $\mathcal{H}^\infty = \mathcal{V} \oplus \mathcal{W}$ one can associate several spectral invariants to this spectral data.
\begin{itemize}
 \item The wave-trace $\mathrm{tr}(e^{ t \mathcal{L}_Z}|_{\ker(\Box)})$ is the distribution in $\mathcal{D}'(\R)$ defined by assigning the test function $\varphi \in C^\infty_0(\R)$ the value 
 $$(\mathrm{tr}(e^{-\rmi  t H}),\varphi) = \sum_{k=1}^N m_k \hat \varphi(\beta_k) + \sum_j \hat\varphi(\lambda_j).$$
 \item The heat trace $\mathrm{tr}(e^{ t \mathcal{L}_Z^2}|_{\ker(\Box)})$ is defined for $t>0$ as
 $$
  \mathrm{tr}(e^{ - t H^2}) = \sum_{k=1}^N m_k e^{-t \beta_k^2} + \sum_j e^{-t \lambda_j^2}.
 $$
 \item The spectral zeta function $\zeta(s)$ is defined as the meromorphic continuation of the function 
 $$
  \zeta(s) = \mathrm{tr}((H^2)^{-s}) =  \sum_{k=1, \beta_k\not=0}^N m_k \beta_k^{-2s} + \sum_j  |\lambda_j|^{-2s},
 $$ 
 for some choice of branch for the logarithm in the expression $\beta_k^{-2s}= e^{-2 s \log(\beta_k)}$.
\end{itemize}
One can also define the counting function $N_c(\lambda) = \#\{ j \mid 0 \leq \lambda_j \leq \lambda\}$. Since we are only interested  in the asymptotic behavior of the counting function as $\lambda \to \infty$ we can disregard the finite dimensional space in the computation of the leading coefficients.

By \cite{SZ18} the wave trace has an isolated singularity at $t=0$ which has an expansion into homogeneous distributions. More precisely, at $t=0$ we have a singularity expansion
$$
  \mathrm{tr}(e^{-\rmi  t H}) \sim  c_0\, \mu_{n-1}(t) + c_1\,  \mu_{n-2}(t) + c_2 \, \mu_{n-3}(t) + \ldots
$$
in the sense that the distribution
$$
 \mathrm{tr}(e^{-\rmi  t H}) - \left( c_0 \,\mu_{n-1}(t) + c_1 \, \mu_{n-2}(t) + c_2 \, \mu_{n-3}(t) + \dots + c_k \, \mu_{n-k-1}(t) \right)
$$
is $C^{k-n+1}$ near $t=0$ for $k>n-1$. In fact, we will see that the above is also of order $O(|t|^{k+1-n})$ as $t \to 0$.
Here the homogeneous distribution $\mu_\alpha(t), \alpha \in \C$ of degree $-\alpha$ is defined by analytic continuation of the holomorphic family of tempered distributions
$$
 \mu_\alpha(t) = \cos(\frac{\pi}{2}\alpha) \Gamma\left(\frac{\alpha+1}{2} \right) |t|^{-\alpha}
 $$
in the half plane $\Re(\alpha)<1$. For $\Re(\alpha)>0$ the $\mu_\alpha$ can be expressed by the oscillatory integrals
$$
 \mu_\alpha(t) = \frac{\Gamma\left(\frac{\alpha+1}{2}\right)}{2 \Gamma(\alpha)}\int_{-\infty}^\infty \mathrm{e}^{-\rmi t \tau} |\tau|^{\alpha-1} \mathrm{d} \tau,
$$
i.e. as Fourier transforms of the $L^1_\mathrm{loc}$ functions $\frac{\Gamma\left(\frac{\alpha+1}{2}\right)}{2 \Gamma(\alpha)} |\tau|^{\alpha-1}$
(we refer to \cite[Vol I]{HoI-IV} for further properties of these distributions).

\begin{rem}
 The distributions $\mu_\alpha$ defined above differ from the homogeneous distributions in \cite{SZ18} by the factor $\frac{\Gamma(\frac{\alpha+1}{2})}{\Gamma(\alpha)}$.
 In contrast to \cite{SZ18} we are here interested also in cases when $\alpha$ is a non-positive integer. With this normalisation the $\mu_\alpha$ is an entire family of homogeneous distributions that is non-vanishing at the integers. Thus, $\mu_{-j}$ spans the one-dimensional space of even homogeneous distributions of degree $j$ for every $j \in \mathbb{Z}$. The $\Gamma$-factors prevent the formation of poles or zeros.
\end{rem}

The coefficient $c_0$ was computed in \cite{SZ18} 
in terms of the symplectic volume $\mathrm{Vol}(\mathcal{N}_{Z \leq 1})$ of the set $\mathcal{N}_{Z \leq 1} = \{\gamma \mid g(\dot \gamma, Z) \leq 1\}$ in the symplectic manifold $\mathcal{N}$ of scaled null geodesics in $T^*M$. Note that $g(\dot \gamma, Z)$ is constant along any null-geodesic.
The first coefficient is given by (see \cite{SZ18}, Theorem 1.3)
$$
c_0 = 2 \pi ^{\frac{1}{2}-n} \Gamma
   \left(\frac{n+1}{2}\right)	 \mathrm{Vol}(\mathcal{N}_{Z \leq 1}) = 2 
  \frac{1}{(2 \pi)^{n-1}} \frac{\Gamma (n)}{\Gamma(\frac{n}{2})}  \mathrm{Vol}(\mathcal{N}_{Z \leq 1}),
$$
which results in the Weyl law (\cite{SZ18}, Cor. 1.5)
$$
 N_c(\lambda) = \frac{1}{(2 \pi)^{n-1}} \mathrm{Vol}(\mathcal{N}_{Z \leq 1}) \lambda^{n-1} + O(\lambda^{n-2} ).
$$
One can express the volume $\mathrm{Vol}(\mathcal{N}_{Z \leq 1})$ in terms local invariant quantities and finds 
(\cite{SZ18}, Equ. (15))
$$
 \mathrm{Vol}(\mathcal{N}_{Z \leq 1}) = \mathrm{Vol}(B_{n-1}) \int_\calK |Z|^{-n} \der \mathrm{Vol}_\calK.
$$

In this paper we show that in principle all the wave-trace invariants can be computed from the Hadamard coefficients of the Hadamard expansion for the Green's functions of $\Box$ near the diagonal. In the first instance this results in an expression for the coefficient $c_2$ in terms of local Lorentzian geometric data.

\begin{theorem} \label{mainth}
 The coefficient $c_1$ vanishes and the coefficient $c_2$ is given by the integral $\int_\Sigma c_2(x) \der \mathrm{Vol}_h(x)$, where
 $c_2(x)= \frac{1}{2}\pi ^{-n/2}  \tilde c_2(x)$ and 
 \begin{align*}
\tilde c_2(x) &=
\left(\frac{1}{6} \mathrm{scal}(x) - W(x) \right) N(x) |Z|^{-n+2} 
- \frac{n-2}{6 }\mathrm{Ric}(Z,Z) N(x) |Z|^{-n}\\+ &\frac{n(n-2)}{12} N(x) g(\nabla_Z Z,\nabla_Z Z) |Z|^{-n-2}+ \frac{1}{3}\mathrm{Ric}(\nu,Z) |Z|^{-n+2} +\frac{n-2}{3} |Z|^{-n} g(\nabla_Z^2 Z,\nu).
\end{align*}
Here $\nu$ denotes the future pointing unit normal vector field along $\Sigma$.\\
The above integral is therefore a spectral invariant of the stationary spacetime.
\end{theorem}

Remarkably the local expression $c_2(x)$ is not a priori invariantly defined on the space of Killing orbits (see Sec. \ref{invariance}) as the last two terms depend on the choice of Cauchy hypersurface. However, we show in Sec. \ref{invariance} that the integral can be rewritten as an integral over locally invariant quantities. The computation is not entirely straightforward.

In case the manifold $(M,g)$ is an ultrastatic spacetime, i.e. when $g = -\der t^2 + h$ is of product form we have
$N=1,\, \nabla_Z Z=0,\, \mathrm{Ric}(Z,Z)=\mathrm{Ric}(\nu,Z)=0$ and therefore the expression for $\tilde c_2$ reduces to the well-known and much simpler second heat kernel coefficient (see for example \cite{MR1032631}, Theorem 1.1)
$$
 \tilde c_2(x) = \left(\frac{1}{6} \mathrm{scal}(x) - W(x) \right),
$$ 
where $\mathrm{scal}(x)$ coincides with the scalar curvature of $(\Sigma,h)$.

As a consequence one therefore obtains that the heat trace $\mathrm{tr}(e^{ t \mathcal{L}_Z^2})$ has an asymptotic expansion of the form
 $$
  \mathrm{tr}(e^{ -tH^2}) = \sum_{k=0}^m a_k t^{-\frac{n-1}{2}+\frac{k}{2}} + o(t^{-\frac{n-1}{2}+\frac{m}{2}} )
 $$
 as $t \to 0_+$. For the coefficients one obtains $a_k = 2^{1 + k - n} \sqrt{\pi} c_k$ and therefore
 $$
  a_0 = \frac{2}{(4 \pi)^{\frac{n-1}{2}}} \int_\calK |Z|^{-n} \der \mathrm{Vol}_\calK, \quad a_1=0, \quad a_2 = \frac{2}{(4 \pi)^{\frac{n-1}{2}}} \int_\Sigma \tilde c_2(x) \der \mathrm{Vol}_h.
 $$
 This reproduces the usual heat kernel coefficients in case the metric is of product type.

Similarly, the spectral zeta function has a meromorphic continuation to the complex plane with at most simple poles as a consequence of the above expansion. This is a direct consequence of
$$
 \sum_j  |\lambda_j|^{-2s} = \frac{1}{\Gamma(s)} \int_0^\infty   t^{s-1} ( \sum_j e^{-\lambda_j^2 t} ) \der t 
 $$
which holds for $\Re(s)$ sufficiently large. 
Since $\sum_j e^{-\lambda_j^2 t}= \mathrm{tr}(e^{- t H^2}) - \sum_{k=1}^N m_k e^{-t \beta_k^2}$ has an asymptotic expansion near $t=0$ this formula provides a meromorphic continuation for the zeta function in the usual manner.

\begin{example}[The rotating round sphere with $N=1$ revisited]
It is illustrative to see what the above theorems simplify to in Example \ref{examplerot}.
In this case the shift vector field $w$ and the Laplace operator on the sphere commute, which makes it possible to explicitly solve the spectral problem.
For the wave operator one obtains
$$
 \Box_g u = (\partial_t-\Omega \partial_\varphi)^2 u - \Delta_h u,
$$
where $\Delta_h$ denotes the Laplace-Beltrami operator of $(S^2,h)$.
The eigenvectors $u_\lambda$ of $Z$ on $\mathrm{ker}(\Box_g)$ corresponding to non-zero eigenvalues are of the form
$$
 u(t,\vartheta,\varphi)= e^{-\rmi \lambda t} Y_{\ell m}(\vartheta,\varphi),
$$
where $Y_{\ell m}$ are the usual spherical harmonics and 
$(\lambda + \Omega m)^2 = \ell(\ell+1)$. There is a non-trivial Jordan block in the decomposition of $Z$ formed by the functions $1$ and $t$ corresponding to the generalised eigenspace with eigenvalue zero.
The rest of the spectrum of $Z$ on $\mathrm{ker}(\Box_g)$ is therefore given by
$$
 \lambda = -\Omega m \pm \sqrt{\ell(\ell+1)}, \quad \ell = 1,2,\ldots, \quad m = -\ell,-\ell+1, \ldots,\ell-1,\ell.
$$
One computes $|Z| =(1-\Omega^2\sin^2\vartheta)^\frac{1}{2}$,
$\nabla_Z Z=-\frac12\,\mathrm{grad}\bigl(g(Z,Z)\bigr)=-\Omega^2\sin\vartheta\cos\vartheta\,\partial_\vartheta$, $\mathrm{Ric}(Z,Z)=\Omega^2\sin^2\vartheta$, and $\mathrm{Ric}(Z,\nu)=0$, $g(\nabla_Z^2 Z, \nu) =0$. 
The local expansion coefficient from Theorem \ref{mainth} is therefore
\begin{align*}
c_2(\vartheta)
&=
\frac{1}{2\pi^{3/2}}
\Bigg[
\frac13(1-\Omega^2\sin^2(\vartheta))^{-1/2}
-\frac{\Omega^2\sin^2(\vartheta)}{6}(1-\Omega^2\sin^2(\vartheta))^{-3/2}\\
&\hspace{4cm}
+\frac{\Omega^4\sin^2(\vartheta)\cos^2(\vartheta)}{4}(1-\Omega^2\sin^2(\vartheta))^{-5/2}
\Bigg].
\end{align*}
Integrating this one obtains
$$
\int_{S^2} c_2\,\der \mathrm{Vol}_h
=
\frac{1}{2\sqrt{\pi}}
\left[
\frac{\mathrm{arctanh}(\Omega)}{\Omega}
+
\frac{1-2\Omega^2}{3(1-\Omega^2)}
\right].
$$
The coefficient $c_0$ computes to 
$\frac{8}{\sqrt{\pi}} \frac{1}{1-\Omega^2}$ and the heat expansion therefore specialises to
\begin{align*}
 \sum_{\pm} \sum_{\ell=0}^\infty \sum_{m=-\ell}^\ell e^{-t (\Omega m \pm \sqrt{\ell(\ell+1)})^2} = \frac{2}{1-\Omega^2} t^{-1} + \frac{1}{2} \left(\frac{1-2
   \Omega ^2}{3 \left(1-\Omega
   ^2\right)}+\frac{\mathrm{arctanh}(\Omega )}{\Omega
   }\right) + O(t)
\end{align*}
as $t \to 0_+$.
In case $\Omega=0$ the function $\frac{\mathrm{arctanh}(\Omega )}{\Omega}$ is understood as a holomorphic function near zero, where it evaluates to $1$.
One then recovers the usual heat expansion on the sphere.
\end{example}

\subsection{Notations and conventions}

The paper uses the following notations and conventions.
The Fourier transform $\hat f$ of $f \in L^1(\R^n)$ will be defined by $$\hat f(\xi) = \int f(x) e^{-\rmi \, x \cdot \xi} \der x$$ where $x \cdot \xi$ is the Euclidean inner product on $\R^n$.
For the metric we choose the sign convention $(-,+,\ldots,+)$ and the d'Alembert operator $\Box$ has principal part in local coordinates $-\sum_{i,k = 1}^n g^{ik} \partial_i \partial_k$ and therefore its principal symbol is
$\sigma_\Box(\xi)=g(\xi,\xi)$. 
Unless otherwise stated functions take values in the complex numbers: For example $C^\infty(M)$ denotes the space of complex valued smooth functions on $M$. The set of real-valued smooth functions on $M$ is denoted by $C^\infty(M,\R)$. For the sake of computations of principal symbols of operators we identify functions on manifolds with half-densities, using the metric and then view the operators as acting between half-densities.

\section{The Hadamard expansion}

In this section we briefly review the Hadamard expansion for the fundamental solution of a normally hyperbolic operator on a globally hyperbolic spacetime $(X,g)$. Constructions are due to  Hadamard \cite{H1,H2} and M. Riesz \cite{R1,R2}. Contemporary expositions include \cite{F,Gun, BK96, Be, BGP}.

By global hyperbolicity there exists a smooth spacelike Cauchy hypersurface $\Sigma \subset X$.
We assume that $\Box: C^\infty(X) \to C^\infty(X)$ is of the form
$\Box = \Box_g + W$, where $- \mathrm{tr}_g \nabla^2$ is the metric d'Alembert operator and $W \in C^\infty(X,\R)$ is a smooth real-valued potential.
Given $s \in \R$ we define $\mathrm{ker}_s(\Box)$
to be the space of solutions with initial data in $H^s_c(\Sigma) \oplus H^{s-1}_c(\Sigma)$. Standard energy estimates imply that this space and its topology do not depend on the choice of Cauchy surface $\Sigma$. 
We will write $\mathrm{ker}(\Box) = \cap_{s>0} \mathrm{ker}_s(\Box)$ for the space of smooth solutions with spacelike compact support.
It is well known that $\Box$ admits unique retarded and advanced fundamental solutions $G_\pm: C^\infty_0(X) \to C^\infty(X)$ and we will denote their difference by
$G = G_+ - G_-$. 
The following exact sequence expresses the main properties of $G$ (see e.g. \cite{BGP}, Theorem 3.4.7)

$$
  \begin{tikzcd}\label{DIAGRAMexact}
  0 \arrow[r] &  C^\infty_c(X) \arrow[r,"\Box"]& C^\infty_c(X) \arrow[r,"G"] & \mathrm{ker}(\Box) \arrow[r,""] &  0.
 \end{tikzcd}
$$

The operator $G$ is a Fourier integral operator of order $-\frac{3}{2}$. This has been shown more generally for certain differences of distinguished parametrices in \cite{DH72}(Equ. (6.5.4) in Theorem 6.5.3) and the retarded and advanced fundamental solutions are such distinguished parametrices (see Section 6.6 of \cite{DH72}).

The singularities of $G$ near the diagonal are described very precisely by the Hadamard expansion that we now explain.
First note that for sufficiently small neighborhood of the diagonal of $X \times X$  the indefinite Lorentzian geodesic square-distance function $\Gamma: \mathcal{U} \to \R$ is well defined as 
$\Gamma(x,y) = -g(\mathrm{exp}_y^{-1}(x),\mathrm{exp}_y^{-1}(x))$.
Let
\begin{equation}
C(\beta,n) := \frac{2^{-n-2\beta}\pi^{\frac{2-n}{2}}}{\Gamma(\beta+\frac{n}{2}) \Gamma(\beta+1)}.
\label{eq:defC}
\end{equation}
The Riesz distributions are then defined for $\Re(\beta)>0$ as the functions 
$$
R^\pm_\beta = 2\,C(\beta,n)\,\Gamma^\beta_\pm,
$$
where $\Gamma^\beta_\pm$ is 
$$
 \Gamma^\beta_\pm(x,y) = \begin{cases} \Gamma^\beta(x,y) & \textrm{if } x \in J^\pm(y),\\ 0 &  \textrm{otherwise},
 \end{cases}
$$
and $J^\pm(y)$ is the causal future/past of the point $y$.
The Riesz distributions admit an analytic continuation as a family of distributions in the parameter $\beta$ to the entire complex plane.
We define $R_\beta = R^+_\beta -  R^-_\beta$. This is an entire family of distributions in $\mathcal{D}'(\mathcal{U})$.
We refer to \cite{BGP} and the Appendix A.3 in \cite{BS} for details.

The Hadamard coefficients are defined in a neighborhood $\mathcal{U}$ of the diagonal of $X \times X$ and form a family $(V_k)_{k \in \N_0}$ of smooth functions $V_k \in C^\infty(\mathcal{U})$. They are uniquely determined by the initial condition
$V_0(x,x) = 1$ and the transport equations
\begin{equation}
 \nabla_{\mathrm{grad}_{(1)} \Gamma} V_k - \left( \tfrac{1}{2}\Box_{0,(1)} \Gamma - n +2k \right) V_k = 2 k \Box_{(1)} V_{k-1}.
\label{eq:transport}
\end{equation}
Here $\Box_{0}$ denotes the scalar wave operator  and the subscript $(1)$ in $\mathrm{grad}_{(1)}$, $\Box_{0,(1)}$, and $\Box_{(1)}$
means that the differential operator acts on the first variable of the subsequent function on $\calU \subset X \times X$. 
The singularity expansion of $G$ on $\mathcal{U}$ is described accurately by the Hadamard expansion. The precise statement is that the $N$-th remainder term of the Hadamard expansion, namely
\begin{align} \label{Hadamard}
  G|_{\mathcal{U}} - \sum\limits_{0 \leq j \leq N+\frac{n-2}{2}} V_{j} \cdot R_{1-\frac{n}2+j}.
\end{align}
is in $C^N(\mathcal{U})$ and is of order $O(\Gamma^N)$ uniformly on compact subsets of $\mathcal{U}$.

The family of Riesz distributions $R^\pm_\beta(x,x_0)$ can therefore also be regarded as the pull-back of the family of Riesz distributions $\mathcal{R}_\beta^\pm$ in Minkowski spacetime under such coordinate systems. We will therefore discuss the Riesz distributions in Minkowski spacetime in more detail.

\subsection{Riesz distributions in Minkowski space}

The Riesz distributions in Minkowski spacetime are of the form,
$$
 R_\beta(x-x') = \mathcal{R}_\beta(x-x'),\quad R^\pm_\beta(x-x') = \mathcal{R}^\pm_\beta(x-x'),
$$
where $\mathcal{R}_\beta$ is a holomorphic family of tempered distributions in Minkowski spacetime. They are constructed by analytic continuation of the homogeneous distributions defined by
$$
  \mathcal{R}^\pm_\beta(x) = \begin{cases} 2\,C(\beta,n)\,\gamma^\beta(x) & \textrm{if } \gamma(x) >0, \pm x_0>0,\\ 0 &  \textrm{otherwise},
 \end{cases}
$$
where $\gamma(x) = x_0^2 -x_1^2 - \ldots -x_{n-1}^2$.
These distributions are holomorphic families of tempered distributions. 
We note that the Fourier transform $\hat{\mathcal{R}}_\beta$ of $\mathcal{R}_\beta$ is an odd distribution which is Lorentz invariant and invariant under spatial reflections. It follows that $\hat{\mathcal{R}}_\beta$ is supported inside the light-cone. This can also be seen from the direct formulae for the Fourier transforms of the Riesz distributions in the literature (e.g.\cite{MR0639526}, p. 85). 
This implies that the family of Riesz distributions $\mathcal{R}_\beta$ is a holomorphic family of distributions in $\mathcal{D}'_\Gamma(\R^n)$, where $\Gamma$ is the closed light cone. This has also been shown independently in great detail in the thesis by Lennart Ronge \cite{Ro23}, Section 4.1. We omit the details.
Note here that $\mathcal{D}'_\Gamma(\R^n)$ denotes the space of distributions with wavefront set contained in $\Gamma$ equipped with a natural topology which is sometimes referred to as the H\"ormander topology.

\subsection{Integration of the Riesz distributions with respect to a timelike vector field}

Suppose that $Z$ is a timelike vector field on $M$. Given $x \in M$ there exists a small relatively compact neighborhood $U$ and an $\epsilon>0$ such that for $t \in (-\epsilon,\epsilon)$ and $x,y \in U$ there exists a unique geodesic connecting $\mathrm{exp}(t Z) x$ and $y$. Hence, the expression
$$
 R_{\beta}(t,x,y) = R_\beta(\exp(t Z) x,y)
$$ 
defines a distribution on $(-\epsilon, \epsilon) \times U \times U$. By the above properties of the Riesz distributions this is a holomorphic family of distributions 
taking values in $\mathcal{D}'_{\tilde \Gamma}((-\epsilon,\epsilon) \times U \times U)$. Here $\tilde \Gamma$ equals
$$
 \{ (t,\tau, x, \xi , y , -\eta) \mid \tau = \xi(Z), (\exp(t Z)^* (x,\xi), y,\eta) 
  \in \Gamma \}.
$$
It follows from wavefront set calculus that for any $f \in C^\infty(-\epsilon,\epsilon)$ the formal expression
$$
 R_{f,\beta}(x,y) = \int_\R R_{\beta}(t,x,y) f(t) \der t
$$
understood as a distributional pairing
defines a smooth function in $C^\infty(U \times U)$ and the map $f \mapsto R_{f,\beta}(x,y)$ is a holomorphic family of distributions in the parameter $\beta$ on $(-\epsilon,\epsilon)$ taking values in $C^\infty(U \times U)$.

We would like to note that in parallel to this work and independently Ronge has obtained such statements in his thesis (\cite{Ro23}, Chapter 4).

\subsection{The Hadamard expansion on a stationary spacetime}

We are now specialising the Hadamard expansion to our setting of a stationary spacetime $(M,g)$.
In this case, purely for symmetry reasons, and by uniqueness of the Hadamard coefficients, we know that the Hadamard coefficients are invariant under the product Killing flow on $e^{t (Z_{(1)}+ Z_{(2)})}= e^{t Z_{(1)}} e^{t Z_{(2)}}$ on $M \times M$. Here, given any differential operator $P$ on $M$ the notation $P_{(1)}$ means the differential operator on $M \times M$ acting on the first variable, and similarly $P_{(2)}$ is acting on the second variable. In terms of (injective) tensor products $C^\infty(M \times M) = C^\infty(M) \otimes_\epsilon C^\infty(M)$ this means for example $Z_{(1)} = Z \otimes 1$ and $Z_{(2)} = 1 \otimes Z$, whereas $Z_{(1)}+ Z_{(2)}= Z \otimes 1 + 1 \otimes Z$.

Note that the Hadamard coefficients are functions defined near the diagonal in $M \times M$. Similarly, the distributions $G_\beta$ are defined near the diagonal only. This neighborhood of the diagonal can be chosen invariant under the product Killing flow.

In the ultrastatic (product) case we have $N=1$ for the Lapse function $N$ and $w=0$ for the shift vector field $w$. In that case the Hadamard coefficients are invariant under time translation in each variable separately. This is however not true in general in the stationary case as we will see below. 
We can write the manifold $M \times M$ as $\R_t \times \Sigma_y \times \R_{t'} \times \Sigma_{y'}$. By symmetry the Hadamard coefficients depend only on the difference $t-t'$. We therefore can write
$$
 V_k(t,y,t',y') = \tilde V_k(t-t',y,y'),
$$
and  $\tilde V_k(t,y,y')$ is defined for $|t|$ sufficiently small and $y,y'$ sufficiently close.
Similarly, the Lorentzian distance function $\Gamma(t,y,t',y')$ is of the form $\tilde\Gamma(t-t',y,y')$
and  $\tilde\Gamma(t,y,y)$, as a distribution in $t$, admits an expansion into homogeneous distributions.
Slightly abusing notations we will denote by $V_k(t,x,y)$ the Hadamard coefficient translated in the $x$-variable by time $t$, in other words $V_k(t,x,y) = V_k(\mathrm{exp}(tZ) x,y)$. This is a well defined smooth function for small $t$ and $x,y$ in a neighborhood of the diagonal in $M$. In case $x,y$ are on a Cauchy surface this coincides with $\tilde V_k(t,x,y)$.
We similarly write $\Gamma(t,x,y)$ for the translated square Lorentzian distance function 
$\Gamma(t,x,y) = \Gamma(\mathrm{exp}(tZ) x,y)$

\section{Computation of the wave-trace}

Following \cite{SZ18} we use the notation $G_s = e^{Z_{(1)} s} G$. We will think of this as a distribution on $\R_s \times M \times M$.
We now recall from \cite{SZ18}, Equ. (25) that the formula for the distributional trace on $\ker(\Box)$ paired with the test function $\varphi \in C^\infty_0(\R)$ is given by
\begin{align*}
   (\mathrm{tr}(e^{Z s} |_{\mathrm{ker}(\Box)}),\varphi) &= \int_\Sigma \int_\R \varphi(s) \left(\partial_{\nu,y} G_s(y,y') - \partial_{\nu,y'}  G_{s}(y,y')\right)|_{y'=y} \der s\;\der V_h(y),
\end{align*}
where $\Sigma$ is a Cauchy surface with future directed normal vector field $\nu$. The integral in $s$ here is understood formally and as a distributional pairing. Not to overload notations we will often write integrals for distributional pairings and we will also write $Q_y(t)$ even when $Q_y$ is a distribution bearing in mind that this expression needs to be understood formally.

The integral in $y$ does not depend on the chosen Cauchy surface. 
For fixed $y \in \Sigma$ one defines the local wave-trace
$Q_y \in \mathcal{D}'(\R)$ as 
$$
 Q_y(\varphi) = \int_\R \varphi(s) \left(\partial_{\nu,y} G_s(y,y') - \partial_{\nu,y'}  G_{s}(y,y')\right)|_{y'=y} \der s
$$
which we abbreviate as
$$
 Q_y(s) =   \left(\partial_{\nu,y} G_s(y,y') - \partial_{\nu,y'}  G_s(y,y')\right)|_{y'=y},
$$
bearing in mind that this is understood as a distribution.
By the above the wave-trace $\mathrm{tr}(e^{Z s} |_{\mathrm{ker}(\Box)})$ equals the integral over the local trace in the sense that
$$
(\mathrm{tr}(e^{Z s} |_{\mathrm{ker}(\Box)}),\varphi) = \int_\Sigma (Q_y, \varphi)\der V_h(y).
$$
We also have
$$
Q_y(s) = \left(\partial_{\nu,y} G_s(y,y')\right)|_{y'=y}  +   \left(\partial_{\nu,y} G_{-s}(y,y')\right)|_{y'=y}.
$$
It will therefore be sufficient to analyze the distribution 
$$
 L_y(s) = \left(\partial_{\nu,y} G_s(y,y')\right)|_{y'=y}.
$$

In what follows we let $I^\alpha, \alpha \in \R$ be the set of distributions in $u \in \mathcal{D}'(\R)$ that have near $0$ a step one polyhomogeneous expansion in terms of homogeneous distributions. Namely, for any $k \in \N$ with $k>\alpha$ there exist  homogeneous distributions $u_{\alpha-j}, j=0,\ldots, k$ of degree $j-\alpha$ on $\R$ with the property that
$$
 u(t) -  \left(  u_{\alpha}(t) +  u_{\alpha-1}(t) +  u_{\alpha-2}(t) + \dots + u_{\alpha-k}(t) \right)
$$
is a continuous function of order $o(t^{k-\alpha})$ in $t$ near zero.
We then write $u \sim v \mod I^\alpha$ if $u-v \in I^\alpha$. If $u$ depends on an additional parameter $y \in \Sigma$ we similarly write $I^m_\Sigma$ for the set of distributions with such expansions as above, where the coefficients $c_j$ depend continuously on $y$ and the order $o(t^{k-\alpha})$ of the remainder term is uniform on compact subsets of $\Sigma$.
In case $u$ is an even distribution in $I^\alpha$ this means there exist numbers $(c_k)_{k \in \N_0}$ with
$$
 u \sim c_0 u_{\alpha}(t) +  c_1 u_{\alpha-1}(t) + c_2 u_{\alpha-2}(t) + \dots + c_k u_{\alpha-k}(t) \quad \mathrm{mod} \quad I^{\alpha-k-1}.
$$

We will compute the first two terms of the polyhomogeneous expansion in $s$ of $Q_y(s)$ using
\begin{equation}
 Q_y(s) = L_y(s) + L_y(-s).
\end{equation}

The expansion of this distribution can conveniently be computed from the Hadamard expansion, which in fact shows that the $y$-dependent family of distributions $L_y(s)$ is in $I^{n-1}_\Sigma$.
This computation can be done in local coordinates and we will therefore fix a coordinate system as follows.
We choose a point $x' \in \Sigma \subset M$ and normal coordinates $(x^0,x^1,\ldots,x^{n-1})$ such that $\partial_{0}|_{x=x'} = \partial_\nu|_{x=x'}$ and the $\partial_j$ are tangent to $\Sigma$ at $x'$.

First the unit-normal derivative $\partial_\nu$ is given by $\frac{1}{N} \left(\partial_t - w^j \partial_j \right)|_{x=x'}$, where we set $w^0=0$.
In other words, we can write $Z|_{x=x'} = \partial_t |_{x=x'}=N \partial_\nu |_{x=x'} + w^j \partial_j|_{x=x'}$.
Further we need some observations about the nature of normal coordinates. 
In these coordinates the metric tensor has the following small $|x|$ expansion
$$
 g_{jk} = \eta_{jk} -\frac{1}{3} R_{j \ell k m}(0) x^\ell x^m + O(|x|^3),
$$
where $\eta_{jk} = \mathrm{diag}(-1,1,\ldots, 1)$ is the Minkowski metric.
The volume distortion function $\mu_{x}$ relative to the point $x'$ is given by $\mu_{x} = |\mathrm{det}(g)|^{\frac{1}{2}}$ when the metric is expressed in normal coordinates about the point $x'$.
This means for small $|x|$ we have
\begin{gather} \label{voldist}
 \mu_x = 1 - \frac{1}{6}\mathrm{Ric}_{k \ell}(0) x^k x^\ell  + O(|x|^3).
\end{gather}
We would now like to express the Killing field $\partial_t$ in terms of the normal coordinates. Observe that at the point $x'$ we have $\partial_t |_{x=x'} = (N(x') \partial_0 + w^j(x') \partial_j)|_{x=x'}$. Let $x(t)$ be the flow line of the Killing flow, i.e. $x(t) = \mathrm{exp}(t Z) x'$. Then the coordinates $x^j(t)$ of this point satisfy
\begin{align*}
 \quad x^j(0) = 0, \quad  &\frac{d x^j(t)}{dt}|_{t=0} =  (N(x'), w(x')), \quad (\frac{d^2 x^j(t)}{dt^2} )|_{t=0}=  
  (\nabla_Z Z(x'))^j,\\
  &(\frac{d^3 x^j(t)}{dt^3} )|_{t=0}=  
  (\nabla_Z \nabla_Z Z(x'))^j.
\end{align*}

We have used in the last step that in normal coordinates the metric Christoffel symbols $\Gamma^\ell_{m n}$
satisfy
$$
 \Gamma^\ell_{m n}(x) = -\frac{1}{3} \eta^{\ell k}\left( R_{k m n j}(0) + R_{k n m j}(0) \right) x^j + O(|x|^2)
$$
as $|x| \to 0$, and $R(\partial_\ell,Z,Z,Z)=0$.
This gives
$$
 x^j(t) =  (N(x'), w(x')) t + \frac{t^2}{2} (\nabla_Z Z)^j + \frac{t^3}{6} (\nabla_Z \nabla_Z Z)^j + O(|t|^4). 
$$
as $t \to 0$.

The first Hadamard coefficient $\tilde V_0(t,y,y')$ is given by
$$
 \tilde V_0(t,y,y') = \mu_{(t,y)}^{-\frac{1}{2}}(0,y').
$$
Here the volume distortion function $\mu_{t,y}$ is given by $\mu_{t,y} = |\mathrm{det}(g)|^\frac{1}{2}$ when the metric is expressed in normal coordinates about the point $x'=(0,y')$.
Using \eqref{voldist} we obtain the following useful formulae in the normal coordinate system
\begin{align}
  V_0(t,y,y) &= 1 + \frac{1}{12} \mathrm{Ric}(Z,Z)(y) t^2 + O(|t|^3) \label{vnullexp}   \\
  \partial_{\nu,y}  V_0(t,y,y')|_{y'=y} &= \frac{1}{6}\mathrm{Ric}(\nu,Z)(y) t + O(|t|^2)    \end{align}
  as $t \to 0$. 
In normal coordinates we also have
$$
 \mathrm{grad}_x \Gamma(x,x') = -2 x^j \partial_{x^j}
$$
which results for small $t$ in
$$
 (\partial_\nu  \Gamma)(t,y,y) = 2 x^0(t) + O(t^4)
$$

$$
 (\partial_\nu  \Gamma)(t,y,y) = 2 N(y) t - g(\nabla_Z Z,\nu)(y) t^2 - \frac{1}{3} g(\nabla_Z^2 Z, \nu)(y) t^3 + O(|t|^4).
$$

As before we write $|Z|^2 = -g(Z,Z)= (N(y)^2 - |w(y)|_h^2)>0$.

By a similar computation we have the small $|t|$ expansion
$$
  \Gamma(t,y,y) = |Z|^2 t^2 + \frac{1}{12} g(\nabla_Z Z,\nabla_Z Z)(y) t^4 + O(|t|^5)  
$$
where we have used that $Z$ is Killing and therefore $g(Z,\nabla_Z Z) =0$. By the antisymmetry properties of the Riemann curvature tensor the terms of the form $R(\partial_0,Z,Z,Z)$ vanish. Moreover, $g(Z,\nabla_Z^2 Z) = -g(\nabla_Z Z,\nabla_Z Z)$.
For $\beta \in \C$ this gives
\begin{align} \label{Gammapower} 
 \Gamma(t,y,y)^\beta = |Z|^{2\beta} |t|^{2 \beta} + \frac{\beta}{12} g(\nabla_Z Z,\nabla_Z Z) |Z|^{2\beta-2} |t|^{2 \beta+2} + O(|t^{2 \beta+3}|).
\end{align}

Furthermore, the diagonal value of the first Hadamard coefficient $V_1(x,x)$ is given by
$$
 V_1(x,x) = \frac{ \mathrm{scal}(x)}{6} - W(x).
$$

\subsection{The leading and the subleading term}

We can now use the expansion \eqref{Hadamard} to obtain the first two homogeneous components of the local wave-trace. 

\begin{align*}
 L_y(s) &= \partial_{\nu_y} \tilde G(s,y,y') |_{y'=y} \\&\sim (\partial_{\nu_y} V_0(s,y,y') |_{y'=y}) R_{1-\frac{n}{2}}(s,y,y) +  V_0(s,y,y) \partial_{\nu_y} R_{1-\frac{n}{2}}(s,y,y') |_{y'=y} \\&+
  V_1(s,y,y) \partial_{\nu_y} R_{2-\frac{n}{2}}(s,y,y') |_{y'=y} \quad \mod I^{n-4}_\Sigma.
\end{align*}
To keep the notation shorter we write $(\partial_\nu R_\beta)(s,x,x)$ for $\partial_{\nu,x} R_\beta(s,x,y) |_{y=x}$.
Then
\begin{align} \label{Lx}
 &L_y(s)
 \sim  \frac{1}{6}\mathrm{Ric}(\nu,Z)(y)\; s\; R_{1-\frac{n}{2}}(s,y,y) +  (\partial_\nu R_{1-\frac{n}{2}})(s,y,y)
 \\
 &+ \frac{1}{12}\mathrm{Ric}(Z,Z) \;s^2\; (\partial_\nu R_{1-\frac{n}{2}})(s,y,y) + (\frac{1}{6} \mathrm{scal}(y) -W(y)) (\partial_{\nu}R_{2-\frac{n}{2}})(s,y,y) \quad \mod I^{n-4}_\Sigma. \nonumber
 \end{align}
We have the following formula for $\partial_\nu R_{\beta+1}(s,y,y) = \frac{1}{4 \beta + 2 n} (\partial_\nu \tilde \Gamma)(s,y,y) R_\beta(s,y,y)$ and therefore
\begin{align*}
 \partial_\nu & R_{\beta+1}(s,y,y) \\ &\sim \frac{1}{2 \beta +  n} R_\beta(s,y,y) \left( N(y) s - \frac{1}{2} g(\nabla_Z Z,\nu)(y) s^2 -\frac{1}{6} g(\nabla^2_Z Z,\nu)(y) \cdot s^3 + \ldots \right),
\end{align*}
where the expansion in the bracket is understood as an asymptotic expansion of step one homogeneous distributions.
Here and in the following the expressions need to be understood as distributions in the variable $s$, weakly meromorphic in $\beta$. For example, $|s|^\beta$ and $\mathrm{sign}(s) |s|^\beta$ are understood as a meromorphic family of distributions by analytic continuation.
We therefore have
\begin{align*}
 (\partial_\nu R_{\beta+1})(s,y,y) &\sim \mathrm{sgn}(s)\frac{2}{2 \beta + n } C(\beta,n)  \Gamma(s,y,y)^\beta \\ &\times \left( N(y) s - \frac{1}{2} g(\nabla_Z Z,\nu)(y) s^2 - \frac{1}{6} g(\nabla^2_Z Z,\nu)(y) \cdot s^3 +\ldots \right).
\end{align*}
Using \eqref{Gammapower} we thus obtain
\begin{align}
 (\partial_\nu R_{\beta+1})(s,y,y) &\sim  \frac{2}{2 \beta + n } C(\beta,n) N(y) |Z|^{2\beta} |s|^{2 \beta+1}  \nonumber\\&-
 \frac{1}{2 \beta + n } C(\beta,n) g(\nabla_Z Z,\nu)(y) |Z|^{2\beta} \mathrm{sgn}(s) |s|^{2 \beta+2} 
 \\&+  \frac{\beta}{6(2 \beta +  n)} C(\beta,n) N(y) g(\nabla_Z Z,\nabla_Z Z) |Z|^{2\beta-2} |s|^{2 \beta+3} \nonumber\\
 &-  \frac{1}{3( 2\beta +  n) } C(\beta,n) g(\nabla^2_Z Z,\nu) |Z|^{2\beta} |s|^{2 \beta+3} \nonumber\\
  &+ \ldots. \nonumber 
\end{align}
As before we use
$$
 \mu_{\beta}(s) = \frac{\Gamma\left(\frac{\beta+1}{2} \right)}{2 \Gamma(\beta)} \int e^{-\rmi s \tau} |\tau|^{\beta-1} \der \tau
$$
for the Fourier transform of $\frac{\Gamma\left(\frac{\beta+1}{2} \right)}{2 \Gamma(\beta)}   |\tau|^{\beta-1}$
and note that as meromorphic families of homogeneous distributions
$$
 \mu_{-\beta}(s) = \cos(\frac{\pi}{2} \beta)\Gamma\left(\frac{-\beta+1}{2} \right) |s|^{\beta}. 
$$
This means
\begin{align*}
 \frac{2}{2 \beta + n } C(\beta,n) |s|^{2 \beta+1} &= -\frac{C(\beta,n)}{\sin(\pi \beta) \Gamma\left(-\beta \right)} \frac{2}{2 \beta + n } \mu_{-2\beta-1}(s),\\
  \frac{1}{2 \beta + n } C(\beta,n) |s|^{2 \beta+3} &= \frac{C(\beta,n)}{\sin(\pi \beta) \Gamma\left(-\beta-1\right)} \frac{1}{2 \beta + n } \mu_{-2\beta-3}(s),
\end{align*}
and in particular, using Euler's reflection formula for the Gamma function, we obtain
\begin{align*}
 \frac{2}{2 \beta + n } C(\beta,n) |s|^{2 \beta+1}|_{\beta \to -\frac{n}{2}} &=\pi ^{-n/2} \mu_{n-1}(s),\\
  \frac{1}{2 \beta + n } C(\beta,n) |s|^{2 \beta+3}|_{\beta \to -\frac{n}{2}}  &= -\pi ^{-n/2} \frac{n-2}{4} \mu_{n-3}(s), \\
  \frac{2}{2 \beta + n } C(\beta,n) |s|^{2 \beta+1}|_{\beta \to 1-\frac{n}{2}} &= \frac{\pi ^{-n/2}}{4 } \mu_{n-3}(s).
\end{align*}

Using this the terms in \eqref{Lx} then expand as follows.

\begin{align*}
  \frac{1}{6} \mathrm{Ric}(\nu,Z) &\,s\,   R_{1-\frac{n}{2}}(s,y,y) \sim  \frac{1}{3}\mathrm{Ric}(\nu,Z) C(\beta,n)  |s|  \Gamma(s,y,y)^\beta |_{\beta=1-\frac{n}{2}} \\&=  \frac{\pi ^{-n/2} }{12 } \mathrm{Ric}(\nu,Z) |Z|^{2-n} \mu_{n-3}(s) \quad \mod I^{n-4}_\Sigma, 
\end{align*}
\begin{align*}
   (\partial_\nu &  R_{1-\frac{n}{2}})(s,y,y) \sim {\pi ^{-n/2} }N(y) |Z|^{-n}  \mu_{n-1}(s) + \tilde c\, \tilde \mu_{n-2}(s)\\
  &+  \frac{\pi ^{-n/2} \frac{n}{2}(\frac{n}{2}-1)}{12} N(y) g(\nabla_Z Z,\nabla_Z Z) |Z|^{-n-2} \mu_{n-3}(s) \\ &
  + \frac{\pi ^{-n/2} (\frac{n}{2}-1)}{6}|Z|^{-n} g(\nabla^2_Z Z,\nu)\mu_{n-3}(s)\quad \mod I^{n-4}_\Sigma, 
\end{align*}
where $\tilde \mu_{n-2}(s)$ is an odd homogeneous distribution and $\tilde c$ is a coefficient that will cancel out in the expansion of $Q_x$ and we therefore do not consider it further.
Another term arises from the expansion \eqref{vnullexp} of $V_0(s,x,x)$ in the $s$-variable:
\begin{align*}
\frac{1}{12}&\mathrm{Ric}(Z,Z) \;s^2\; (\partial_\nu  R_{1-\frac{n}{2}})(s,y,y) \\&\sim -\frac{1}{6} \,\mathrm{Ric}(Z,Z)  \frac{\pi ^{-n/2}  \left(\frac{n}{2}-1\right)}{2} N(y)  |Z|^{-n}  \mu_{n-3}(s) \quad \mod I^{n-4}_\Sigma\\
& \sim -\frac{n-2}{24} \,\mathrm{Ric}(Z,Z)  {\pi ^{-n/2}} N(y)  |Z|^{-n}  \mu_{n-3}(s) \quad \mod I^{n-4}_\Sigma.
\end{align*}

Finally,
\begin{align*}
 &\left (\frac{1}{6} \mathrm{scal}(y) - W(y) \right)(\partial_{\nu} R_{2-\frac{n}{2}})(s,y,y) \\&\sim
 \frac{1}{4}\left(\frac{1}{6} \mathrm{scal}(y) - W(y)\right) N(y)  |Z|^{-n+2} {\pi ^{-n/2} } \mu_{n-3}(s) \quad \mod I^{n-4}_\Sigma.
 \end{align*}

Collecting all the terms into homogeneous ones we get
\begin{align*}
  Q_x(s) \sim c_0 (x) \mu_{n-1}(s) + c_2(x) \mu_{n-3}(s) \quad \mod I^{n-4}_\Sigma.
  \end{align*}

The leading term is therefore
\begin{align}
  c_0(x)&= 2  {\pi ^{-n/2}}N(x) |Z|^{-n}  \\&= 2 
  \frac{1}{(2 \pi)^{n-1}}\mathrm{Vol}(\mathrm{B}_{n-1}) \frac{\Gamma (n)}{\Gamma(\frac{n}{2})} N(x) |Z|^{-n}  . \nonumber
\end{align}
This is in agreement with the formula obtained in \cite{SZ18}, Equ.~(15) by a different method. Here $\mathrm{B}_{n-1}$ is the volume of the unit ball in $\R^{n-1}$.
\begin{align*}
  c_2(x)&=   \frac{\pi ^{-n/2}}{2} \times \tilde c_2(x),\\
\tilde c_2(x) &=
\left(\frac{1}{6} \mathrm{scal}(x) - W(x) \right) N(x) |Z|^{-n+2} \\
&- \frac{n-2}{6}\mathrm{Ric}(Z,Z) N(x) |Z|^{-n} +  \frac{n(n-2)}{12} N(x) g(\nabla_Z Z,\nabla_Z Z) |Z|^{-n-2}\\&+ \frac{1}{3}\mathrm{Ric}(\nu,Z) |Z|^{-n+2} +\frac{n-2}{3} |Z|^{-n} g(\nabla_Z^2 Z,\nu). 
\end{align*}

This proves Theorem \ref{mainth}.

\section{Wave-trace coefficients in terms of local invariant quantities} \label{invariance}

If $\Box$ is the usual d'Alembert operator the spectrum of $Z$ on $\ker \Box$ has several natural invariance properties. A straightforward one is the natural transformation property under Killing-diffeomorphisms. This generalizes the fact that the spectrum of the Laplacian on isometric Riemannian manifolds is the same.
A particular class of such diffeomorphism are the local-gauge transformations of $M$ regarded as a principal $\R$-bundle.

Given a function $f: M \to \R$ that is invariant under the Killing flow it is easy to see that the family of diffeomorphisms $\Psi_\epsilon: M \to M$ given by
$$
 \Psi_\epsilon = \mathrm{exp}(\epsilon f Z)
$$
is equivariant in the sense that it commutes with the Killing flow. This diffeomorphism is a local gauge transformation if one thinks of $M$ is an $\R$-principal bundle over the space of Killing orbits.
If $M = \R \times \Sigma$ and the metric has standard stationary form the diffeomorphism is of course simply given as
$$
 \Psi_\epsilon(t,x) = (t + \epsilon f(x), x) 
$$
We will assume this in this section without loss of generality.
We can use the pulled back metric
$g_\epsilon = \Psi_\epsilon^* g$ to construct the d'Alembert operator $\Box_\epsilon = \Box_{g_\epsilon}$.
It is then clear that the spectrum of $Z$ on $\ker \Box_\epsilon$ does not depend on $\epsilon$.
This implies that the spectral invariants must also be invariant under this transformation.
Given $f$ it is straightforward to compute the infinitesimal variations with respect to the above group of diffeomorphism at zero. For a quantity $F_\epsilon$ we write $\delta F = \partial_\epsilon F_\epsilon |_{\epsilon=0}$ and
 obtain in local coordinates
\begin{align*}
  \delta h_{jk} &=  h_{\ell k} w^\ell \frac{\partial f}{\partial x^j} + h_{\ell j} w^\ell \frac{\partial f}{\partial x^k},\\
  \delta \sqrt{\mathrm{det}(h)} &= w^k \frac{\partial f}{\partial x^k}\sqrt{\mathrm{det}(h)},\\
  \delta w_j &=  -(N^2-w^2)  \frac{\partial f}{\partial x^j},\\
   \delta N &= -N w^k \frac{\partial f}{\partial x^k},\\
   \delta R_{00}&=0,  \delta R_{k0} = R_{00} \frac{\partial f}{\partial x^k},\\
   \delta \partial_\nu &= \frac{1}{N} w^k \frac{\partial f}{\partial x^k} \partial_t +  N h^{jk} \frac{\partial f}{\partial x^j} \partial_k
  \end{align*}
This shows that the following quantities are invariant,
\begin{align*}
    \delta (N^2 - |w|_h^2)&=0,\\
    \delta \left(N \sqrt{\mathrm{det}(h)}\right)=0.
\end{align*}
In terms of $u = |Z|$ and $\theta = \der t - u^{-2} h_{jk} w^k \der x^j$ we have
$$
 \delta u =0,\quad \delta \theta = \der f.
$$

As explained in the introduction these transformation properties can also be obtained directly from the invariant description of $M$ as a principal bundle. 
The invariance of $N \sqrt{\mathrm{det}(h)}$ follows from the fact that this is the invariant volume form on the space of Killing orbits.

If we look at the local coefficient in terms of invariant quantities we see that in
\begin{align*}
\tilde c_2(x) &=
\left(\frac{1}{6} \mathrm{scal}(x) -W(x) \right) N(x) |Z|^{-n+2} \\
&- \frac{n-2}{6}\mathrm{Ric}(Z,Z) N(x) |Z|^{-n} + \frac{n(n-2)}{12} N(x) g(\nabla_Z Z,\nabla_Z Z) |Z|^{-n-2}\\&+ \frac{1}{3} |Z|^{-n}\left( \mathrm{Ric}(\nu,Z) |Z|^2 +(n-2) g(\nabla_Z^2 Z,\nu) \right)
\end{align*}
all but the last term are local invariants when multiplied by the Riemannian volume form $\der \mathrm{Vol}_h$.
The last is not a priori invariant. Since the integral of this expression over $\Sigma$ is a spectral invariant, it is invariant under the above transformation, or in other words must be invariant under a change of Cauchy surface. It must therefore be possible to re-write this expression in terms of invariant quantities. 

To see that this is indeed the case we recall the Pohozaev-Schoen identity (Prop 1.4 in \cite{MR929283}). It states that on a Riemannian manifold $(M,g)$ with boundary $\partial M$ and a conformal Killing field $X$ we have
$$
 \frac{n-2}{2n}\int (\mathcal{L}_X R) \,\der \mathrm{Vol}_{g} = \int_{\partial M} (\mathrm{Ric}(X,\nu) - \frac{1}{n} R \cdot g(X,\nu)) \der \sigma_g.
$$
The proof of this identity is based on a local computation, using the second Bianchi and the divergence theorem. It therefore carries over without change to Lorentzian manifolds.
We will apply this identity to the metric $\tilde g = |Z|^{-2} g$. The induced metric on the Cauchy surface is $\tilde h = |Z|^{-2} h$ and therefore $\der \sigma_{\tilde g} = |Z|^{-n+1} \der \sigma_{g}$. Since $Z$ is a Killing field for $\tilde g$ the above identity implies that
\begin{align*}
 I_M &= \int_\Sigma (\tilde{\mathrm{Ric}}(Z,\tilde \nu) - \frac{1}{n} \tilde R \cdot \tilde g(Z,\tilde \nu)) \der \sigma_{\tilde g} \\&=
  \int_\Sigma (\tilde{\mathrm{Ric}}(Z,\nu) - \frac{1}{n} \tilde R \cdot  \tilde g(Z, \nu)) |Z|^{-n+2}  \der \sigma_{g} 
\end{align*}
does not depend on the Cauchy surface $\Sigma$ and is therefore an invariant of the stationary spacetime $(M,g,Z)$.
Here we have used that $g(\nu,\nu) = \tilde g (\tilde \nu, \tilde \nu) =|Z|^{-2} g(\tilde \nu, \tilde \nu)$ and therefore
$\tilde \nu =|Z| \nu$.
By the formula for conformal change of the Ricci curvature (\cite{MR2371700}, Theorem 1.159 (d)) we have
$$
 \tilde{\mathrm{Ric}}(Z,\nu) =\mathrm{Ric}(Z, \nu) + \Box \varphi g(Z,\nu) - (n-2) g^{-1}(\der \varphi, \der \varphi) g(Z,\nu)
 -(n-2) (\mathrm{Hess}\, \varphi)(Z,\nu)
$$
where $\varphi = - \log |Z|$.
Since $g(Z,\nu) = -N$ we obtain
\begin{align*}
 &\int_\Sigma \mathrm{Ric}(\nu,Z) |Z|^{-n+2} \der \mathrm{Vol}_h \\&=
 I_M + \int_\Sigma \left(\Box \varphi - (n-2) g^{-1}(\der \varphi, \der \varphi) \right) |Z|^{-n+2}  N \der \mathrm{Vol}_h \\&-\int_\Sigma
 \frac{1}{n}  \tilde R \cdot  |Z|^{-n}  N \der \mathrm{Vol}_h \\&+\int_\Sigma (n-2) (\mathrm{Hess}\, \varphi)(Z,\nu)   |Z|^{-n+2}   \der \mathrm{Vol}_h  
\end{align*}
The conformally transformed scalar curvature is given by
\begin{align*}
 \tilde R = R  (N(x)^2 - |w(x)|_h^2) &+2 (n-1) |Z|^2 \Box \varphi -(n-2)(n-1)  |Z|^2 g^{-1}(\der \varphi,\der\varphi). 
\end{align*}
Hence, collecting all terms one has
\begin{align*}
 \int_\Sigma \mathrm{Ric}(\nu,Z)  |Z|^{-n+2} &\der \mathrm{Vol}_h- \int_\Sigma (n-2) (\mathrm{Hess}\,\varphi)(Z,\nu)   |Z|^{-n+2}   \der \mathrm{Vol}_h  \\&=
 I_M + \frac{n-2}{n}\int_\Sigma \left(-\Box \varphi - g^{-1}(\der \varphi, \der \varphi) \right)  |Z|^{-n+2}   N \der \mathrm{Vol}_h \\&-\int_\Sigma
 \frac{1}{n}  R \cdot   |Z|^{-n+2}  N \der \mathrm{Vol}_h  \\&=
 I_M + \frac{n-2}{n}\int_\Sigma \left(\Box |Z| \right)  |Z|^{-n+1} N \der \mathrm{Vol}_h \\&-\int_\Sigma
 \frac{1}{n}  R \cdot   |Z|^{-n+2}  N \der \mathrm{Vol}_h 
\end{align*}

Since $Z \varphi=0$ the Hessian is given by
$$
 (\mathrm{Hess}\, \varphi)(Z,\nu) = -(\nabla_\nu Z) \varphi = -\frac{1}{2}  |Z|^{-2} (\nabla_\nu Z) g(Z,Z),
$$
bearing in mind that $g(Z,Z) = -|Z|^2$.
Since $Z$ is a Killing field we have  for any vector fields $X,Y$ the Killing equation $g(\nabla_X Z, Y) + g(X,\nabla_Y Z)=0$.
This shows 
$$
 (\nabla_\nu Z) g(Z,Z) = 2 g(\nabla_{(\nabla_\nu Z)},Z) = -2 g(\nabla_\nu Z, \nabla_Z Z) = 2 g(\nu, \nabla_{(\nabla_Z Z)} Z).
$$
Since $[Z,\nabla_Z Z]=0$ and by torsion-freeness of the Levi-Civita connection this finally gives the equation
$$
 (\mathrm{Hess}\, \varphi)(Z,\nu) |Z|^2 = -g (\nu, \nabla_Z^2 Z).
$$
Hence,
\begin{align*}
  \int_\Sigma  |Z|^{-n}\left( \mathrm{Ric}(\nu,Z) |Z|^2 +(n-2) g(\nabla_Z^2 Z,\nu) \right) \der \mathrm{Vol}_h  \\= 
  I_M +\int_\Sigma \left( \frac{n-2}{n} |Z|^{-1}\Box |Z|  - \frac{1}{n}  R  \right)  |Z|^{-n+2}  N \der \mathrm{Vol}_h,
\end{align*}
which is invariant as claimed.

\section*{Acknowledgement}

The results of this paper were part of an ongoing research collaboration between the two authors to generalize aspects of Riemannian spectral geometry to settings applicable to general relativity. Most of this paper has been written in collaboration between the first author and the late Steve Zelditch. The first author would like to use this opportunity to express his deep gratitude for this long term collaboration and inspiration.

\nocite{}

\end{document}